\documentclass[11pt]{article}
\usepackage{color}
\usepackage{amsmath,amsthm,amssymb}
\newtheorem{theorem}{Theorem}
\newtheorem{lemma}{Lemma}
\newtheorem{coro}{Corollary}
\newtheorem{defi}{Definition}

\newcommand{\ba}{\begin{array}}
\newcommand{\ea}{\end{array}}
\newcommand{\bt}{\begin{tabular}}
\newcommand{\et}{\end{tabular}}
\newcommand{\btb}{\begin{table}}
\newcommand{\etb}{\end{table}}
\newcommand{\bc}{\begin{center}}
\newcommand{\ec}{\end{center}}
\newcommand{\bea}{\begin{eqnarray}}
\newcommand{\eea}{\end{eqnarray}}
\newcommand{\Bea}{\begin{eqnarray*}}
\newcommand{\Eea}{\end{eqnarray*}}
\newcommand{\beq}{\begin{equation}}
\newcommand{\eeq}{\end{equation}}
\newcommand{\Beq}{\begin{equation*}}
\newcommand{\Eeq}{\end{equation*}}

\begin{document}
\title{ \bf A class of Solvable Lie algebras}
\author{Yan Wang,\quad Ran Cui,\quad ShaoQiang Deng\thanks{Corresponding author.
{\em E-mail address:} dengsq@nankai.edu.cn.}\\
{\footnotesize\em School of Mathematical
Sciences and LPMC, }\\
{\footnotesize\em Nankai University, Tianjin 300071, China
}}
\date{}
\maketitle

\baselineskip 18pt

\noindent {\bf Abstract:}
\medskip
All finite-dimensional indecomposable solvable Lie algebras $\frak
g$, having the filiform Lie algebra $Q_{2m+1}$ as the nilradical,
are studied and classified. It turns out that the dimension of
$\frak g$ is at most $\dim{Q_{2m+1}}+2$.


\section{Introduction}

\qquad Solvable Lie algebras play an important role in Lie theory. A
fundamental theorem due to Levi ([1]) asserts that any finite
dimensional Lie algebra $\frak g$ can be decomposed into a
semidirect sum
\begin{equation*}
\frak g={\frak s}\dotplus{\frak r},\quad [{\frak s}, {\frak
s}]={\frak s},\quad [{\frak r}, {\frak r}]={\frak r},\quad [{\frak
s}, {\frak r}]\subseteq \frak r,
\end{equation*}
where $\frak s$ is semisimple and $\frak r$ is the radical of $\frak
g$, i.e. its maximal solvable ideal. Semisimple Lie algebras over
the field of complex numbers ${\bf C}$ have been completely
classified by Cartan ([2]), over the field of real numbers ${\bf R}$
by Gantmacher ([3]). So, to complete the classification of all Lie
algebras amounts to classify the solvable
Lie algebras. However, for solvable Lie algebras such
classifications exist only for low dimensions (e.g. [4], [5]). It
seems to be very difficult to proceed by dimension in the
classification of Lie algebra $\frak g$ beyond $\dim{\frak g}=6$, in
particular solvable ones. It is however possible to proceed by
structure.

The purpose of this paper is to classify a certain type of finite
dimensional solvable Lie algebras over the field ${\bf C}$. It is
well known that any solvable Lie algebra $\frak g$ has a uniquely
defined nilradical $N$, i.e. maximal nilpotent ideal, satisfying
 $\dim N\geq \frac{1}{2}\dim \frak g.$
Hence we can consider a given nilpotent Lie algebra of dimension $n$
as the nilradical and then find all of its extensions to solvable
Lie algebras. Here, we regard filiform Lie algebra $Q_n$ as the
nilradical of the solvable Lie algebra. When $n=2m$, $Q_{2m}$ is a
kind of naturally graded nilpotent Lie algebras, which has been
discussed in [6]. So we only need to concentrate on the case
$N=Q_{2m+1}\ (m\geq 2)$ with Lie brackets
\begin{eqnarray}
&&[e_1, e_k]=e_{k+1},\quad 2\leq k\leq 2m, \label{1}\\
&&[e_k, e_{2m+2-k}]=(-1)^k e_{2m+1},\quad 2\leq k\leq m,\nonumber
\end{eqnarray}
over the basis $\{e_1, e_2, \cdots, e_{2m+1}\}$. In previous
articles the classification has been performed on the following nilpotent Lie
algebras: Heisenberg algebras $H_n(n\geq1)$ ([7]),  Abelian Lie
algebras $a_n(n\geq1)$ ([8]), `triangular' Lie algebras $t_n(n\geq
2)$ ([9]), filiform Lie algebras $L_n(n\geq 4)$ ([10]),
quasifiliform  algebras $q_n(n\geq 4)$ ([11]).

The paper is organized as follows. In Section 2, we recall some basic
concepts related to Lie algebras and the definition of filiform algebras.
In Section 3, we calculate the outer derivations and automorphisms for
$Q_{2m+1}$, respectively, and obtain their matrix forms.
Finally, in Section 5, we present the classification of the solvable
Lie algebras with nilradical $Q_{2m+1}$.

\section{Preliminaries}
\subsection{Basic concepts}

\qquad Let $\frak g$ be a Lie algebra over a field $\bf F$. The
derived series ${\frak g}={\frak g}^{(0)}\supseteq {\frak g}^{(1)}
\supseteq \cdots \supseteq {\frak g}^{(k)}\supseteq \cdots$ is
defined as
$$
{\frak g}^{(0)}={\frak g},\quad {\frak g}^{(k)}=[{\frak g}^{(k-1)},
{\frak g}^{(k-1)}].
$$
If there exists $k\in \bf N$ such that ${\frak g}^{(k)}=0$, then $\frak g$
is called a solvable Lie algebra.

The lower central series ${\frak g}={\frak g}^0\supseteq {\frak g}^1
\supseteq \cdots {\frak g}^k \supseteq \cdots$ is defined as
\begin{equation*}
{\frak g}^0={\frak g},\quad {\frak g}^k=[{\frak g}, {\frak g}^{k-1}].
\end{equation*}
If there exists $k\in \bf N$ such that ${\frak g}^k=0$, then $\frak
g$ is called a nilpotent Lie algebra and the smallest integer
satisfying ${\frak g}^k=0$ is called the nilindex of $\frak g$.
\begin{defi}$^{[12]}$
The Lie algebra $\frak g$ is called filiform if ${\rm dim}{\frak
g}^k=n-k-1$ for $k\geq 1$.
\end{defi}

A derivation $D$ of a Lie algebra $\frak g$ is a linear transformation of $\frak g$  such that
\begin{equation}\label{2}
D([x,\ y])=[D(x), y]+[x, D(y)].
\end{equation}
If there exists an element $z\in \frak g$ such that
\begin{equation*}
D={\rm ad}z,\quad i.e., \quad D(x)=[z, x],\quad \forall x\in \frak g,
\end{equation*}
then the derivation is called an inner derivation. A derivation
which is not inner will be called an outer derivation.

\subsection{Basic classification methods}

\qquad Any solvable Lie algebra $\frak g$ contains a unique
nilradical $N$. We will assume that $N$ is given, that is, we know
the Lie brackets of $N$ with respect to the basis $\{e_1, e_2,
\cdots, e_{n}\}$:
\begin{equation}\label{3}
[e_a, e_b]=\sum\limits_{c=1}^{n}{C^c_{ab}e_c}.
\end{equation}
We study the problem of the classification of all the solvable Lie
algebras $\frak g$ with nilradical $N$. This can be achieved by
adding further elements $h_1, h_2, \cdots, h_t$ to the basis $\{e_1,
e_2, \cdots, e_n\}$ to form a basis of $\frak g$. Since the derived
algebra of such a solvable Lie algebra is contained in the
nilradical ([13]), i.e.
\begin{equation}\label{4}
[\frak g, \frak g]\subseteq  N,
\end{equation}
we have
\begin{eqnarray}
&&[h_i, e_a]=\sum\limits_{b=1}^n(H_i)^b_ae_b, \quad 1\leq i\leq t,\  1\leq a\leq n,\label{5}\\
&&[h_i, h_j]=\sum\limits_{a=1}^{n}r^a_{ij}e_a, \quad 1\leq i, j\leq
t. \label{6}
\end{eqnarray}
Obviously, each $h_i$ is correspondent to one matrix $H_i$.

The matrix elements of the matrices $H_i$ must satisfy certain linear relations following
from the Jacobi identities between the elements $(h_i, e_a, e_b)$. In addition, since $N$ is
the maximal nilpotent ideal of $\frak g$, the matrices $H_i$ must satisfy another condition:
no nontrivial linear combination of them is a nilpotent matrix, that is, they are linearly
nil-independent.

Now, we will consider the matrices $H_i$ from the other hand. Let us
consider the adjoint representation of $\frak g$, restrict it to the
nilradical $N$ and find ${\rm ad}|_N(h_i)$. It follows from the
Jacobi identities that ${\rm ad}|_N(h_i)$ is a derivation of $N$.
Moreover, it is an outer derivation. In other words, finding all
sets of matrices $H_i$ in (\ref{5}) satisfying the Jacobi identities
is equivalent to finding all sets of outer nil-independent
derivations of $N$:
\begin{equation*}
D_1={\rm ad}|_N(h_1), \cdots, D_t={\rm ad}|_N(h_t).
\end{equation*}
Furthermore, in view of (\ref{4}), the commutators $[D_i, D_j]$ must
be inner derivations of $N$. This requirement determines the
structure constants $r^a_{ij}$.

Different sets of derivations may correspond to isomorphic Lie algebras. In terms of
(\ref{5}) and (\ref{6}), it means that a classification of the solvable Lie algebra
$\frak g$ thus amounts to a classification of the matrices $H_i$
and the structure constants $r^a_{ij}$ under the following transformations:
\begin{eqnarray}
\tilde{h}_i &=& \sum\limits_{j=1}^{t}G_{ij}h_j + \sum\limits_{a=1}^nS_{ia}e_a, \label{7} \\
\tilde{e}_a &=& \sum\limits_{b=1}^nR_{ab}e_b, \label{8}
\end{eqnarray}
where $G$ is an invertible $t\times t$ matrix, $S$ is a $t\times n$
matrix and the invertible $n\times n$ matrix $R$ must be  chosen so
that the Lie brackets (\ref{3}) are preserved.

\section{The outer derivations of $Q_{2m+1}$}
\qquad The nilpotent Lie algebra $Q_{2m+1}$ is defined by the Lie
brackets (\ref{1}) in the introduction. From Section 2.2, we know
that determining the set of outer derivations of $N$ plays an
important role in classification. Now, we will describe explicitly the
set of outer derivations of $Q_{2m+1}$.

First, we introduce some notations. Let
$$[e_i, e_j]=B_{ij}e^\alpha,$$
where $B_{1j}=(\overbrace{0, \cdots, 0}^{j}, 1, 0, \cdots, 0)\
(2\leq j\leq 2m),\ B_{i,2m+2-i}=(0, \cdots, 0, (-1)^i)\ (2\leq i\leq
m)$, other ${B_{ij}}^\prime$s are zero row vectors and
$e^\alpha=(e_1, e_2, \cdots, e_{2m+1})^T$. Let $D$ be an
outer derivation of $Q_{2m+1}$. Suppose
$D(e_i)=\sum\limits_{k=1}^{2m+1}d_{ik}e_k$. We have
$D=(d_{ij})_{i,j=1}^{2m+1}$. Apply condition (2) to basis elements
$e_i$ and $e_j$, that is
$$D[e_i,e_j]=[De_i,e_j]+[e_i, De_j].$$
We  obtain
\begin{eqnarray}\label{9}
B_{ij}D=d_{i1}B_{1j}+\cdots+d_{i,2m+1}B_{2m+1,j}+d_{j1}B_{i1}+\cdots+d_{j,2m+1}B_{i,2m+1}.
\end{eqnarray}
In the following, we discuss (\ref{9}) with given $i$ and $j$ in
order to get the matrix $D$.

(1)$\ i=1,\ 2\leq j\leq m: (d_{j+1,1}, \cdots, d_{j+1,2m+1})=(0, 0,
d_{j2}, \cdots, d_{11}+d_{jj}, \cdots,
d_{j,2m}+(-1)^{j+1}d_{1,2m+2-j})$;

(2)$\ i=1,\ m+2\leq j\leq 2m-1: (d_{j+1,1}, \cdots,
d_{j+1,2m+1})=(0, 0, d_{j2}, \cdots, d_{11}+d_{jj}, \cdots,
d_{j,2m}+(-1)^jd_{1,2m+2-j})$;

(3)$\ i=1,\ j=m+1: (d_{m+2,1}, \cdots, d_{m+2,2m+1})=(0, 0,
d_{m+1,2}, \cdots, d_{11}+d_{m+1,m+1}, \cdots, d_{m+2,2m})$.\\
From the above three cases, we get
\begin{align*}
&d_{j,k}=0, \qquad\qquad\qquad\quad         k<j\ (3\leq j\leq 2m),\\
&d_{j,k}=d_{2,k-j+2}, \qquad\qquad    (3\leq j<k\leq 2m),\\
&d_{11}+d_{jj}=d_{j+1,j+1}\qquad (2\leq j\leq 2m-1).
\end{align*}
Somewhat differently, when $2\leq j\leq m$, we have
$d_{j+1,2m+1}=d_{j,2m}-(-1)^jd_{1,2m+2-j}$ and when $m+2\leq j\leq
2m-1$, we have $d_{j+1,2m+1}=d_{j,2m}+(-1)^jd_{1,2m+2-j}$.

(4)$\ i=1,\ j=2m: (d_{2m+1,1}, \cdots, d_{2m+1,2m+1})=(0, \cdots, 0,
d_{11}+d_{12}+d_{2m,2m})$;

(5)$\ i+j=2m+2,\ i\leq m: ((-1)^id_{2m+1,1}, \cdots,
(-1)^id_{2m+1,2m+1})=(0,
\cdots, 0, (-1)^id_{ii}+(-1)^id_{jj})$;\\
From  (4) and (5), we have $d_{2m+1,k}=0\ (1\leq k\leq 2m),\
d_{2m+1,2m+1}=d_{11}+d_{12}+d_{2m,2m}$ and
$d_{ii}+d_{jj}=d_{2m+1,2m+1}$. Let $d_{11}=\alpha$ and
$d_{22}=\beta$. It is easy to get $d_{kk}=(k-2)\alpha+\beta\ (3\leq
k\leq 2m),\ d_{2m+1,2m+1}=(2m-2)\alpha+2\beta$ and
$d_{12}=\beta-\alpha$.

(6)$\ i=2,\ j=3: (0, \cdots, 0)=(0, 0, 0, d_{21}, 0, \cdots, 0,
d_{2,2m-1}+d_{3,2m})$. \\
Obviously, $d_{21}=0$. Up to now, we can get that $D$ is an upper
triangular matrix in the basis $\{e_1, e_2, \cdots, e_{2m+1}\}$ of
$Q_{2m+1}$. A further computation shows that

$d_{2s}=0\ (m+2\leq s\leq 2m-1$ and $s$ is odd $)$ when $i=2$ and
$j=2m+2-s$;

$d_{2t}=0\ (3\leq t\leq m+1)$ when $i=m+3-t$ and $j=m+1$.\\
Therefore we have $d_{j+1,2m+1}=(-1)^jd_{1,2m+2-j}\ (m+2\leq j\leq
2m-1)$ and $d_{m+2,2m+1}=d_{2,m-1}=0$.

Finally, according to the equivalence of classification (\ref{7}),
we can simplify the matrices of the outer derivations by changing
the basis in the space span$\{D_1, \cdots, D_t\}$. Denote
$\tilde{h}_k=h_k+\sum\limits_{i=1}^{2m+1}a_ie_i\ (1\leq k\leq t)$.
We have the following lemma.

\begin{lemma}
The outer derivation of $Q_{2m+1}$ is upper triangular matrix in the
basis $\{e_1, e_2, \cdots, \\e_{2m+1}\}$. The diagonal elements
satisfy:
$$ d_{kk}=(k-2)\alpha+\beta\ (3\leq k\leq 2m),\
d_{2m+1,2m+1}=(2m-2)\alpha+2\beta,$$ \\
where $d_{11}=\alpha$ and
$d_{22}=\beta$. The specific form of the matrix $D$ is
\begin{equation}\label{10}
\left( \arraycolsep=0.5pt
\begin{array}{cccccccccccc}
\alpha &\beta-\alpha  &0             &\cdots    &d_{1,m+1} &0 &0         &0         &\cdots  &0          &0                  &d_{1,2m+1}-d_{2,2m+1}\\
0      &\beta         &0             &\cdots    &0         &0 &d_{2,m+3} &0         &\cdots  &0          &d_{2,2m}           &0 \\
0      & 0            &\alpha+\beta  &\cdots    &0         &0 &0         &d_{2,m+3} &\cdots  &d_{2,2m-2} &0                  &d_{2,2m}\\
       &              &              &\ddots    &          &  &          &          &\ddots  &           &d_{2,2m-2}         &0\\
       &              &              &          &\ddots    &  &          &          &        &\ddots     &                   &d_{2,2m-2}\\
       &              &              &          &          &  &          &          &        &           &           &\vdots\\
       &              &              &          &          &  &          &          &        &           &                   &d_{2,m+3}\\
       &              &              &          &          &  &          &          &        &           &                  &d_{1,m+2}\\
       &              &              &          &          &  &          &          &\ddots  &           &                   &0\\
       &              &              &          &          &  &          &          &        &\ddots     &                   &\vdots\\
0      &0             &0             &\cdots    &0         &0 &0         &0         &\cdots  &0          &(2m-2)\alpha+\beta &0\\
0      &0             &0             &\cdots    &0         &0 &0         &0         &\cdots  &0          &0                  &(2m-2)\alpha+2\beta
\end{array}
\right)
\end{equation}
when $m$ is odd and
\begin{equation}\label{11}
\left(
\arraycolsep=0.8pt
\begin{array}{ccccccccccc}
\alpha &\beta-\alpha  &0             &\cdots    &d_{1,m+1} &0         &0         &\cdots  &0          &0                  &d_{1,2m+1}-d_{2,2m+1}\\
0      &\beta         &0             &\cdots    &0         &d_{2,m+2} &0         &\cdots  &0          &d_{2,2m}           &0 \\
0      & 0            &\alpha+\beta  &\cdots    &0         &0         &d_{2,m+2} &\cdots  &d_{2,2m-2} &0                  &d_{2,2m}\\
       &              &              &\ddots    &          &          &          &\ddots  &           &d_{2,2m-2}         &0\\
       &              &              &          &\ddots    &          &          &        &\ddots     &                   &d_{2,2m-2}\\
       &              &              &          &          &          &          &        &           &           &\vdots\\
       &              &              &          &          &          &          &        &           &                   &d_{2,m+2}-d_{1,m+2}\\
       &              &              &          &          &          &          &\ddots  &           &                   &0\\
       &              &              &          &          &          &          &        &\ddots     &                   &\vdots\\
0      &0             &0             &\cdots    &0         &0         &0         &\cdots  &0          &(2m-2)\alpha+\beta &0\\
0      &0             &0             &\cdots    &0         &0         &0         &\cdots  &0          &0                  &(2m-2)\alpha+2\beta
\end{array}
\right)
\end{equation}
when $m$ is even.
\end{lemma}

\section{The  automorphism  of $Q_{2m+1}$ }

\qquad As we explained in Section 2.2, to find all solvable Lie
algebras with nilradical $Q_{2m+1}$, we must find all nonequivalent
nil-independent outer derivation sets $\{D_1, \cdots, D_t\}$ of
$Q_{2m+1}$. The equivalence is determined by the transformations
(\ref{7}) and (\ref{8}). The transformation (\ref{7}) has been
considered in Section 3. Here, we consider the transformation
(\ref{8}). Denote the transformation by $\sigma$, we have
$\tilde{e}_a=\sigma(e_a)$. It should be noted that $\sigma$ must
remain the Lie brackets, that is
\begin{equation}\label{12}
[\sigma(e_i), \sigma(e_j)]=\sigma[e_i, e_j].
\end{equation}
Obviously, $\sigma$ is an automorphism of $Q_{2m+1}$. In this section,
 we will  determine the form of automorphism of $Q_{2m+1}$.

Let $\sigma(e_i)=\sum\limits_{k=1}^{2m+1}a_{ik}e_k$. Then
$$
\left(
\begin{array}{c}
\sigma e_1\\
\sigma e_2\\
\vdots\\
\sigma e_{2m+1}
\end{array}
\right)
=
\left(
\begin{array}{cccc}
a_{11}     &a_{12}     &\cdots &a_{1,2m+1}\\
a_{21}     &a_{22}     &\cdots &a_{2,2m+1}\\
\vdots     &\vdots     &       &\vdots\\
a_{2m+1,1} &a_{2m+2,2} &\cdots &a_{2m+1,2m+1}
\end{array}
\right)
\left(
\begin{array}{c}
e_1\\
e_2\\
\vdots\\
e_{2m+1}
\end{array}
\right)
$$
In short, $(\sigma e_1, \cdots, \sigma e_{2m+1})^T=Ae^\alpha$. By
equation (\ref{12}), we have
$$[a_{i1}e_1+\cdots+a_{i,2m+1}e_{2m+1},\ a_{j1}e_1\\
+\cdots+a_{j,2m+1}e_{2m+1}]=B_{ij}Ae^\alpha.$$ Furthermore,
\begin{align}
&(0, 0, a_{i1}a_{j2}-a_{i2}a_{j1}, \cdots,
a_{i1}a_{j,2m-1}-a_{i,2m-1}a_{j1},
(a_{i1}a_{j,2m}-a_{i,2m}a_{j1})+\nonumber\\
&(a_{i2}a_{j,2m}-a_{i,2m}a_{j2})-\cdots+(-1)^m(a_{im}a_{j,m+2}-a_{i,m+2}a_{jm}))=B_{ij}A.\label{13}
\end{align}
Similarly to the method of discussing the derivation $D$, we discuss
the equation (\ref{13}) with given $i$ and $j$. It is easy to get
that $A$ is an upper triangular matrix and
\begin{align*}
&a_{12}=q-p,\ a_{kk}=p^{k-2}q\ (3\leq k\leq 2m),\\
&a_{2m+1,2m+1}=p^{2m-2}q^2,\ a_{ij}=p^{i-2}a_{2,j-i+2}\ (3\leq
i<j\leq 2m),
\end{align*}
where $p=a_{11}$ and $q=a_{22}$. In addition, we can get the
following three results:

(1)$\ a_{2s}=0\ (x<m+2$ and $s$ is odd $)$ when $i=m+3-s$ and $j=m+1$;

(2)$\ a_{2s}=0\ (x\leq m+1$ and $s$ is even $)$ when $i=2$ and $j=2m+2-s$;

(3)$\ a_{2s}=0\ (m+2\leq x\leq 2m-1$ and $s$ is odd $)$ when $i=3$ and $j=2m+1-s$.\\
Now, we can give the simplest form of the matrix A,
$$\left(
\arraycolsep=1.5pt
\begin{array}{cccccccccccc}
p   &q-p  &a_{13} &\cdots    &a_{1,m+2} &a_{1,m+3} &a_{1,m+4}  &\cdots   &a_{1,2m-1}  &a_{1,2m}      &a_{1,2m+1}\\
0   &q    &0      &\cdots    &0         &a_{2,m+3} &0          &\cdots   &0           &a_{2,2m}      &a_{2,2m+1} \\
0   &0    &pq     &\cdots    &0         &0         &pa_{2,m+3} &\cdots   &pa_{2,2m-2} &0             &a_{3,2m+1}\\
    &     &       &\ddots    &          &          &           &         &            &p^2a_{2,2m-2} &a_{4,2m+1}\\
    &     &       &          &\ddots    &          &           &         &            &              &\vdots\\
    &     &       &          &          &          &           &         &            &              &a_{m+1,2m+1}\\
    &     &       &          &          &          &           &         &            &              &0\\
    &     &       &          &          &          &           &\ddots   &            &              &a_{m+3,2m+1}\\
    &     &       &          &          &          &           &         &\ddots      &              &\vdots\\
0   &0    &0      &\cdots    &0         &0         &0          &\cdots   &0           &p^{2m-2}q     &a_{2m,2m+1}\\
0   &0    &0      &\cdots    &0         &0         &0          &\cdots   &0           &0             &p^{2m-2}q^2
\end{array}
\right)
$$
when $m$ is odd and
$$\left(
\arraycolsep=1.5pt
\begin{array}{cccccccccccc}
p   &q-p  &a_{13} &\cdots    &a_{1,m+1} &a_{1,m+2} &a_{1,m+3}  &\cdots   &a_{1,2m-1}  &a_{1,2m}      &a_{1,2m+1}\\
0   &q    &0      &\cdots    &0         &a_{2,m+2} &0          &\cdots   &0           &a_{2,2m}      &a_{2,2m+1} \\
0   &0    &pq     &\cdots    &0         &0         &pa_{2,m+2} &\cdots   &pa_{2,2m-2} &0             &a_{3,2m+1}\\
    &     &       &\ddots    &          &          &           &         &            &p^2a_{2,2m-2} &a_{4,2m+1}\\
    &     &       &          &\ddots    &          &           &         &            &              &\vdots\\
    &     &       &          &          &          &           &         &            &              &a_{m+1,2m+1}\\
    &     &       &          &          &          &           &         &            &              &0\\
    &     &       &          &          &          &           &\ddots   &            &              &a_{m+3,2m+1}\\
    &     &       &          &          &          &           &         &\ddots      &              &\vdots\\
0   &0    &0      &\cdots    &0         &0         &0          &\cdots   &0           &p^{2m-2}q     &a_{2m,2m+1}\\
0   &0    &0      &\cdots    &0         &0         &0          &\cdots   &0           &0             &p^{2m-2}q^2
\end{array}
\right)
$$
when $m$ is even. It should be noted that
$$a_{i,2m+1}=(-1)^{i+1}p^{i-3}qa_{1,2m+3-i}\ (m+2<i<2m+1).$$

Suppose that an automorphism $\sigma$ of $Q_{2m+1}$ makes a transformation
for the basis $\{e_1, e_2,\\
\cdots, e_{2m+1}\}$. Then $\sigma$ can be treated as a
transformation directly on $D$, i.e. $\sigma$ acts on $D$ as $\sigma
D\sigma^{-1}$. Obviously, $D$ is equivalent to $\sigma
D\sigma^{-1}$. Hence, we can further simplify $D$ by using
automorphisms with the above form.

\section{Classification of solvable Lie algebras with the nilradical $Q_{2m+1}$}

\qquad Let $D_i$ be one of the outer derivations of $Q_{2m+1}$. From
Section 3, we know that the diagonal elements of the matrix $D_i$
can be completely determined by $\alpha_i$ and $\beta_i$. Now, we
assume that $\{D_1, D_2, D_3\}$ are linearly independent outer
derivations. According to the theory of linear algebra, the
equations
$$\left\{\begin{array}{l}
 \textrm{$X_1\alpha_1+X_2\alpha_2+X_3\alpha_3=0$}\\
 \textrm{$X_1\beta_1+X_2\beta_2+X_3\beta_3=0$}\\
  \end{array}\right.
$$
has nonzero solutions. Hence $X_1D_1+X_2D_2+X_3D_3$ can be a
nilpotent matrix. In other words, $\{D_1, D_2, D_3\}$ are not
linearly nil-independent. So the number of nil-independent elements
$t$ can be at most two.
\begin{theorem}
Any solvable Lie algebra $\frak g$ with nilradical $Q_{2m+1}$ will
have dimension $\dim{\frak g}=2m+2$ or $\dim{\frak g}=2m+3 $.
\end{theorem}
In the theorems, `solvable' will always mean solvable,
indecomposable, non-nilpotent.

For simplicity, in the matrix $D$, denote $d_{2,k}$ by $d_k\
(m+2\leq k\leq 2m)$, $d_{1,m+1}$ by $a$, $d_{1,2m+1}-d_{2,2m+1}$ by
$b$, $d_{1,m+2}$ by $c\ (\mbox{when}\ m\ \mbox{is odd})$ and
$d_{2,m+2}-d_{1,m+2}$ by $c\ (\mbox{when}\ m\ \mbox{is even})$. In
the following, we will determine the solvable Lie algebra $\frak g$.

\subsection{The case $\dim \frak g=2m+2$}

\qquad The entire structure of $\frak g$ is determined by the matrix
$D$. The Lie brackets of non-nilpotent element $h$ with the basis
$\{e_1,e_2, \cdots, e_{2m+1}\}$ of $Q_{2m+1}$ are given by
$$[h, e_i]=D(e_i)=\sum\limits_{k=1}^{2m+1}d_{ik}e_k.$$
We shall divide our discussion into two cases according to the
values of the parameters $\alpha$ and $\beta$, at least one of which
must be nonzero. Note that, each of the following step is
independent. Only the values of the $(1,2m+1)$-element and $c\ (m$
is even $)$ will be changed when we take other changes. For
simplicity, we still denote them by $b$ and $c$.

(1)$\ \alpha\neq 0$. We rescale $D$ to put $\alpha=1$ and change the
basis of $Q_{2m+1}$
$$\tilde{e}_1=e_1,\ \tilde{e}_i=e_i-\frac{1}{l-2}d_le_{i+l-2}\ (2\leq i\leq 2m+3-l),\
\tilde{e}_j=e_j\ (2m+4-l\leq j\leq 2m+1).$$ Then we can change
$d_{m+3}, d_{m+5}, \cdots, d_{2m}\ (m$ is odd $)$ or $d_{m+2},
d_{m+4}, \cdots, d_{2m}\ (m$ is even $)$ to zero sequentially. At
the same time, the $(1,l)$-element turns into
$\frac{\beta-1}{l-2}d_l$ and the $(2m+3-l,2m+1)$-element turns into
$\frac{1-\beta}{l-2}d_l$, where $l\in\{m+3, m+5, \cdots, 2m\}$ when
$m$ is odd and $l\in\{m+4, m+6, \cdots, 2m\}$ when $m$ is even. If
$l=m+2$, only the $(1,m+2)$-element equals
$\frac{\beta-1}{m}d_{m+2}$. The $(m+1,2m+1)$-element $c$ is changed
but can not be changed to its negative. We still denote it by $c$.

If $\beta\neq2-m$, we continue to change the basis
$\tilde{e}_1=e_1-\frac{a}{(m-2)+\beta}e_{m+1}$, which turns $a$ to
zero.

If $\beta\neq \frac{3-2m}{2}$, a further change of basis
$\tilde{e}_1=e_1-\frac{b}{(2m-3)+2\beta}e_{2m+1}$ turns $b$ to
zero.

Now, all matrix elements in $D$ are zeroes except the diagonal
elements, $c$, the $(1,l)$-element and the $(2m+3-l,2m+1)$-element,
where $l\in\{m+3, m+5, \cdots, 2m\}$ when $m$ is odd and $l\in\{m+2,
m+4, \cdots, 2m\}$ when $m$ is even. It is easily seen that $c$ can not be
changed to zero by any automorphism of $Q_{2m+1}$. In addition, if
$\beta=1$, it is easy to see that the $(1,l)$-element and the
$(2m+3-l,2m+1)$-element are zeroes. If $\beta\neq 1$, let
$$\tilde{e}_1=e_1-\frac{d_l}{(l-3)+\beta}e_l,\quad \tilde{e}_{2m+3-l}=e_{2m+3-l}
+\frac{d_l}{(l-3)+\beta}e_{2m+1}.$$ This change turns the
$(1,l)$-element and the $(2m+3-l,2m+1)$-element into zeroes, where
$l\in\{m+3, m+5, \cdots, 2m\}$ when $m$ is odd and $l\in\{m+4, m+6,
\cdots, 2m\}$ when $m$ is even. But if $l=m+2$, we can only change
the $(1,m+2)$-element to zero. Moreover, if the value of $\beta$
fails to change the corresponding matrix elements to zero, then we can
choose the automorphism
$$\sigma^{\prime}=\rm {diag}(p, p, \cdots, p^{2m})$$
with appropriate $p$ to change them to 1 (except this diagonal
elements).

(2)$\ \alpha=0, \beta\neq 0$. We rescale $D$ to put $\beta=1$ and take the change
$$\tilde{e}_1=e_1-ae_{m+1},$$
which transforms $a$ to zero. We further change the basis
$$\tilde{e}_1=e_1-\frac{b}{2}e_{2m+1},$$
which turns $b$ to zero. If now $D$ is diagonal, it can not be
further simplified. Let us assume that $D$ is not diagonal. All we
can do is to use the transformation $\sigma^{\prime}$ to change one
of the nonzero elements to 1. In summarizing, we have proved:
\begin{theorem}
Seven types of solvable Lie algebras with dimension $2m+2$
exist for $m\geq 2$. They are mutually non-isomorphic and can be represented
as the following:

(1)$\ {\frak g}_{2m+2,1}: [h, e_1]=e_1+(\beta-1)e_2,\ [h, e_i]=(i-2+\beta)e_i\
(2\leq i\leq 2m, i\neq m+1),\ [h, e_{m+1}]=(m-1+\beta)e_{m+1}+\mu e_{2m+1},\
[h, e_{2m+1}]=(2m+2\beta-2)e_{2m+1}.$

(2)$\ \ {\frak g}_{2m+2,2}: [h, e_1]=e_1-me_2,\ [h, e_i]=(i-m-1)e_i\
(2\leq i\leq 2m, i\neq m+1),\ [h, e_{m+1}]=\mu e_{2m+1}$ when $m$ is odd;

$[h, e_1]=e_1-me_2+\nu e_{m+2},\ [h, e_i]=(i-m-1)e_i\ (2\leq i\leq
2m, i\neq m+1),\ [h, e_{m+1}]=\mu e_{2m+1}$ when $m$ is even.

(3)$\ {\frak g}_{2m+2,3}: [h, e_1]=e_1+(1-m)e_2+\mu e_{m+1},\ [h, e_i]=(i-m)e_i\
(2\leq i\leq 2m, i\neq m+1),\ [h, e_{m+1}]=e_{m+1}+\nu e_{2m+1},\
[h, e_{2m+1}]=2e_{2m+1}.$

(4)$\ {\frak g}_{2m+2,4}: [h, e_1]=e_1+(\frac{1}{2}-m)e_2+\mu e_{2m+1},\
[h, e_i]=(i-m-\frac{1}{2})e_i\ (2\leq i\leq 2m, i\neq m+1),\
[h, e_{m+1}]=\frac{1}{2}e_{m+1}+\nu e_{2m+1},\
[h, e_{2m+1}]=e_{2m+1}.$

(5)$\ {\frak g}_{2m+2,5}: [h, e_1]=e_1-e_2,\ [h, e_i]=(i-2)e_i\
(3\leq i\leq 2m, i\neq m+1),\ [h, e_{m+1}]=(m-1)e_{m+1}+\mu e_{2m+1},\
[h, e_{2m+1}]=(2m-2)e_{2m+1}.$

(6)$\ {\frak g}_{2m+2,6}: [h, e_1]=e_1+(\beta-1)e_2-d_{3-\beta}e_{3-\beta},\
[h, e_i]=(i-2+\beta)e_i\ (2\leq i\leq 2m, i\neq m+1, 2m+\beta),\
[h, e_{m+1}]=(m-1+\beta)e_{m+1}+\mu e_{2m+1},\
[h, e_{2m+\beta}]=(2m+2\beta-2)e_{2m+\beta}+d_{3-\beta}e_{2m+1},\
[h, e_{2m+1}]=(2m+2\beta-2)e_{2m+1},$
where $\beta\in\{-2m+3, \cdots, -m-2, -m\}$ when $m$ is odd and
$\beta\in\{-2m+3, \cdots, -m-3, -m-1\}$ when $m$ is even.

(7)$\ {\frak g}_{2m+2,7}: [h, e_1]=e_2,\
[h, e_i]=e_i+\sum\limits_{k=[\frac{m+1}{2}]+1}^m d_{2k}e_{i+2k-2}\ (2\leq i\leq 2m, i\neq m+1),\
[h, e_{m+1}]=e_{m+1}+\mu e_{2m+1},\
[h, e_{2m+1}]=2e_{2m+1}.$\\
In all classification, $\mu=0$ or $1$ and $\nu\in \bf F.$
\end{theorem}

\begin{coro}
${\frak g}_{2m+2,3}$ can not exist for $m=2$.
\end{coro}

\subsection{The case $\dim \frak g=2m+3$}

\qquad Since any nonzero linear combinations of $D^1, D^2$ can not
be nilpotent matrix, there are only one case: $\alpha_1=1,
\beta_1=0$ and $\alpha_2=0, \beta_2=0$.

First, we take $D_1$ to its simplest form by changing the basis of $Q_{2m+1}$:
$$D_1=\left(
\begin{array}{cccccccc}
1                &-1        & 0          &\cdots         &0        & 0       & 0      &0 \\
0                & 0        & 0          &\cdots         &0        & 0       & 0      &0\\
0                & 0        & 1          &\cdots         &0        & 0       & 0      & 0 \\
                 &          &            &\ddots         &         &         &        &\vdots\\
                 &          &            &               &\ddots   &         &        &\mu\\
                 &          &            &               &         &\ddots   &        &\vdots\\
                 &          &            &               &         &         &\ddots  &0\\
0                &0         &0           &\cdots         &0        &0        &0       &2m-2
\end{array}
\right),
$$
where the value of $u$ can only be $0$ or $1$. Let $D_2$ be the form
of matrix (\ref{10}) or (\ref{11}). Since $[D_1, D_2]\in
B=$span$\{\mbox{ad}e_1, \mbox{ad}e_2, \cdots, \mbox{ad}e_{2m+1}\}$,
we compare $[D_1, D_2]$ with matrix $B$ and immediately get
$$D_2=\left(
\begin{array}{cccccccc}
0                &1         & 0          &\cdots         &0        & 0       & 0      &0 \\
0                &1         & 0          &\cdots         &0        & 0       & 0      &0\\
0                &0         & 1          &\cdots         &0        & 0       & 0      & 0 \\
                 &          &            &\ddots         &         &         &        &\vdots\\
                 &          &            &               &\ddots   &         &        &\frac{\mu}{m-1}\\
                 &          &            &               &         &\ddots   &        &\vdots\\
                 &          &            &               &         &         &\ddots  &0\\
0                &0         &0           &\cdots         &0        &0        &0       &2
\end{array}
\right).
$$
In addition, $[D_1, D_2]=0$.

Now, the corresponding Lie bracket of non-nilpotent elements $h_1$
and $h_2$ with $\{e_1, e_2, \cdots, \\e_{2m+1}\}$ can be determined.
It remains to fix the Lie bracket between $h_1$ and $h_2$. Note that
$D_1$ and $D_2$ are commuting matrices, which means the Lie bracket
of $h_1$ and $h_2$ must in the center of $Q_{2m+1}$, i.e. $[h_1,
h_2]= \gamma e_{2m+1}$. The transformation
$$\tilde{h}_1=h_1+\frac{\gamma}{2}e_{2m+1},\ \tilde{h}_2=h_2$$
turns $\gamma$ to zero while leaving other Lie brackets invariant.
In summarizing, we get:

\begin{theorem}
There is only one type of solvable Lie algebra of ${\rm dim}\frak
g=2m+3$ with nilradical $Q_{2m+1}$. Its Lie brackets are given by:
\begin{eqnarray*}
&&[h_1, e_1]=e_1-e_2,\ [h_1, e_i]=(i-2)e_i\ (3\leq i\leq 2m, i\neq m+1),\\
&&[h_1, e_{m+1}]=(m-1)e_{m+1}+\mu e_{2m+1},\ [h_1, e_{2m+1}]=(2m-2)e_{2m+1},\\
&&[h_2, e_1]=e_2,\  [h_2, e_i]=e_i\ (2\leq i\leq 2m, i\neq m+1),\\
&&[h_2, e_{m+1}]=e_{m+1}+\frac{\mu}{m-1}e_{2m+1},\ [h_2, e_{2m+1}]=2e_{2m+1},\\
&&[h_1, h_2]=0,
\end{eqnarray*}
where $u=0$ or $1$.
\end{theorem}

\vskip 6mm
\leftline{\bf REFERENCES}
\vskip 4mm \footnotesize \baselineskip 11pt
\begin{description}

\bibitem[1]]
E. E. Levi, Sula structtura dei gruppi finiti e continui, Atti
Accad. Sci. Torino CI. Sci. Fis. Mat. Natur. 40 (1905) 551-565.

\bibitem[2]]
E. Cartan,  Sur la structure des groups de transformations finis et
continus (Paris: thesis, Nony), in: Oeuvres Completes, Partie I.
Tome 1, 1894, pp. 137-287.

\bibitem[3]]
F. Gantmacher, Rec. Math. (Mat. Sbornik) N. S. 5 (1939) 217-250.

\bibitem[4]]
J. Patera, H. Zassenhaus, Solvable Lie algebras of dimension$\leq 4$
over perfect fields, Linear algebra Appl. 142 (1990) 1-17.

\bibitem[5]]
W. A. de Graaf, Classification of solvable Lie algebras,
Experimental Mathematics 14 (2005) 15-25.

\bibitem[6]]
J. M. Ancochea, R. Campoamor-Stursberg, L. Garcia Vergnolle,
Solvable Lie algebras with naturally graded nilradicals and their
invariants, J. Phys. A: Math. Gen. 39 (2006) 1339-1355.

\bibitem[7]]
J. Rubin, P. Winternitz, Solvable Lie algebras with Heisenberg
ideals, J. phys. A: Math. Gen. 26 (1993) 1123-1138.

\bibitem[8]]
J. C. Ndogmo, P.  Winternitz, Solvable Lie algebras with Abelian
nilradicals, J. phys. A: Math. Gen. 27 (1994) 405-423.

\bibitem[9]]
S. Trembly, P. Winternitz, Solvable Lie algebras with
triangular nilradicals,
J. phys. A: Math. Gen. 31 (1998) 789-806.

\bibitem[10]]
L. $\rm \check{S}$nobl, P. Winternitz, A class of solvable Lie
algebras and their Casimir invariants, J. phys. A: Math. Gen. 38
(2005) 2687-2700.

\bibitem[11]]
Y. Wang, J. Lin, S. Q. Deng, Solvable Lie algebras with
quasifiliform nilradicals, Comm. Algebra, 36 (2008), 1-16.

\bibitem[12]]
M. Goze, Y. Khakimdjanov, Nilpotent Lie algebras, Kluwer Academic Press, 1996.

\bibitem[13]]
N. Jacobson, Lie algebras, Courier Dover Publications, New York, 1979.
\end{description}
\end{document}